# APPLYING NUMEROSITY TO SURREAL INTEGRATION

ILYA SHNITMAN

ABSTRACT

We present a novel framework for measuring the size of discrete subsets of $\mathbb{R}$ using surreal-valued numerosity, which strictly satisfies Euclid's principle that "the whole is greater than a part". By mapping numerosities to surreal numbers via the canonical embedding of Hardy fields, we provide an explicit formula for the full numerosity of sequences, including their infinite, finite and infinitesimal parts, with examples. Then (interpreting germs at infinity as Laplace transforms) we introduce a surreal-valued function that shares many properties with the Dirac Delta distribution and employ it to derive some integration properties of surreal values. Finally, we provide a formula for integrating a surreal-valued function over surreal domain, employing numerosity in a way, similar to Lebesgue measure and derive some novel formulas, connecting surreal numbers via integration. The suggested here method of surreal integration is different from methods suggested by other authors in dropping the property of linearity regarding infinite factors. In our opinion, the struggle to keep this rule intact is the reason for many failed attempts to define the surreal integration. For the suggested formulas, code in Wolfram Mathematica language is provided.

## 1. Mapping numerosities to surreal numbers.

Several authors[1][2][3][4][5] have introduced the concept of a measure of set's size, striving to satisfy the Euclid's principle "the whole is greater than a part" (often calling it "numerosity"), which can be formally represented by the following equation: $N(A \cup B) = N(A) + N(B) - N(A \cap B)$, where $N$ is the proposed measure. The approach described in this paper makes an attempt to do so in the most simple way, yet allowing to explicitly derive values of numerosities in closed form.

Our aim is to define a measure on subsets of $\mathbb{R}$ that would be

- Surreal-valued.
- Equal to cardinality at finite sets.
- Finitely-additive.
- For a bounded set, invariant under translations.
- Finite part is equal to the regularization of the sum over indicator function

*Comb-form of the numerosity.* We formally represent the numerosity of a discrete set $S \subset \mathbb{R}$ with elements $s_n$ as a (generally, divergent) integral $N(S) = \displaystyle\int_{-\infty}^{\infty} \sum_{n=0}^{\infty} \delta(x - s_n)dx$, where $\delta(x)$ is the Dirac Delta distribution. These divergent integrals should be considered as formal symbols until we establish a mapping from them to surreal numbers. Indeed, naturally partial sums of this integral increase by one each time they encounter an element of $S$. We will call integrals of this form the "comb-form" of the numerosity. We assume that numerosity is invariant regarding changing of the sign of the elements of a set, so $\displaystyle\int_{-\infty}^{\infty} \sum_{n=0}^{\infty} \delta(x - s_n)dx = \int_{-1}^{\infty} \sum_{n=0}^{\infty} \delta(x - |s_n|)dx$

*Smooth-form of the numerosity.* Let $\phi_n(x)$ to be a bounded real function with finite support $[a_n, b_n]$, area under the curve $\displaystyle\int_{a_n}^{b_n} \phi_n(x)dx = 1$ and center of mass at $s_n$: $\displaystyle\int_{a_n}^{b_n} x\phi_n(x)dx = s_n$. We postulate that these two integrals are equal: $\displaystyle\int_{-\infty}^{\infty} \sum_{n=0}^{\infty} \delta(x - s_n)dx = \int_{0}^{\infty} \sum_{n=0}^{\infty} (\phi_n(x) + \phi_n(-x))dx$. The later one is one of the smooth-forms of $N(S)$. As with the comb form, smooth form, if diverges, should be considered a formal symbol for now, but we are going to establish mapping between them and surreal numbers.

The surreal numbers can be considered an H-field with natural embedding of the germs of the fuctions (satisfying the Hardy field requirements) at infinity into surreals, such that the germ of identity function $[x]$ is mapped to $\omega$ [6]. In such setting, if $u \in \mathbb{J}$, a purely-infinite surreal number, corresponding to the germ at infinity of the function $f(x)$, we postulate equivalence of the following divergent integral to $u$:

$$f(a) + \int_{a}^{\infty} f'(t)dt = u \tag{1}$$

where $a$ is an arbitrary real number.

In this way, smooth forms of numerosities can be mapped to surreal numbers.

EXAMPLE. Find a subset of $\mathbb{R}$, whose numerosity is $\ln\omega$. The number $\ln\omega$ is purely-infinite, because its Conway normal form $\omega^{1/\omega}$ has no finite and infinitesimal terms. Under canonical embedding and the above-mentioned equivalence,

$$\ln\omega = [\ln x] = \ln 1 + \int_{1}^{\infty} (\ln t)'\, dt = \int_{1}^{\infty} \frac{1}{t}\, dt.$$

Now, we take the last integral and separate it into pieces of area equal to 1, in our case:

$$\int_{1}^{\infty} \frac{1}{t}\, dt = \int_{1}^{e} \frac{1}{t}\, dt + \int_{e}^{e^2} \frac{1}{t}\, dt + \cdots = \sum_{k=0}^{\infty} \int_{e^k}^{e^{k+1}} \frac{1}{t}\, dt.$$

The centers of mass of each integral are located at $s_k = e^{k+1} - e^k$, so the set $\{e^{k+1} - e^k\} : k \in \mathbb{N}_0$ has the numerosity $\ln\omega$.

EXPLICIT FORMULA. The above-mentioned definitions can be expressed in compact form as follows:

$$N(\{a_k\}) = \left(D\Delta^{-1}a_k\right)^{[-1]}(\omega), \tag{2}$$

where $f^{[-1]}$ is the inverse function operator. The formula works not only for purely-infinite numbers but for those containing finite and infinitesimal parts, so we suggest calling the result of applying this formula the full numerosity and after trimming the infinitesimal part, the refined numerosity $\overline{N}(\{a_k\})$. Differences in refined numerosities between different sequences, calculated according to the method suggested in this article coincide with those obtained with the method suggested by James Propp in[7].

MATHEMATICA CODE. Below is the Mathematica code for finding the surreal expression of the numerosity of a non-negative, growing sequence $a_k$ with no accumulation points:

```
a[k_] := k^2;
SolveValues[D[Sum[a[k], k], k] == \[Omega], k] /. C[1] -> 0 //
   FullSimplify // Expand
```

The opposite process, finding a sequence with given numerosity as a surreal expression $S$ can be done with the following code:

```
S := Log[\[Omega]];
DifferenceDelta[Integrate[Normal[SolveValues[S == k, \[Omega]]], k], k] /. C[1]
-> 0 // FullSimplify // Expand
```

In the obtained results only the real positive elements should be considered. The obtained list of numerosities of subsets of reals can be seen in Appendix 1.

## 2. Surreal Delta function

Suppose we want to obtain a functional to count the number of roots of a function on an interval. For a differentiable function $f(x)$, the following expression does the thing, where applicable: $\displaystyle\int_a^b \delta(f(x))|f'(x)|dx$. From now on we will call the operator $\delta_S$, defined as $\delta_S(f(x)) = \delta(f(x))|f'(x)|$ where $x$ is the integration variable, “surreal delta function”. We will further generalize it.

First of all, we should notice that $\delta_S$ behaves just like $\delta(x)$ in all contexts where its argument is the integration variable itself, and not some non-trivial function of it. Particularly, its Laplace and Fourier transforms are the same as of conventional Delta distribution.

Second thing to notice is that it satisfies the following property: $\delta_S(ax) = \delta_S(x)$, which coincides with the similar property of the unit impulse function, unlike the property of the conventional Delta distribution $\delta(ax) = \dfrac{1}{|a|}\delta(x)$.

Let us consider the Fourier transform of the constant function $f(x) \equiv 1$:

$$F(s) = \mathcal{F}[1](s) = \int_{-\infty}^{\infty} e^{-st}dt$$

Formally, the right-hand expression is equal to $0$ at any $s \neq 0$ (by Cesaro regularization). But at $s = 0$ we get the divergent integral $\displaystyle\int_{-\infty}^{\infty} 1dt$. This divergent integral corresponds to the surreal number $2\omega$.

Consequently, we can formally write:

$$F(x) = \mathcal{F}[1](x) = \begin{cases} 2\omega, & \text{if } x = 0 \\ 0, & \text{if } x \neq 0 \end{cases}$$

Applying inverse Fourier transform, $\mathcal{F}^{-1}[F](t) = 1 = 2\pi\mathcal{F}^{-1}[\delta](t) = 2\pi\mathcal{F}^{-1}[\delta_S](t)$ This allows us to generalize the Surreal Delta function this way:

$$\delta_S(x) = \begin{cases} \frac{\omega}{\pi}, & \text{if } x = 0 \\ 0, & \text{if } x \neq 0 \end{cases} \tag{3}$$

We can see that the property $\delta_S(ax) = \delta_S(x)$ holds with this piece-wise definition. This function $\delta_S(x)$ is a full-fledged surreal-valued function, and can be evaluated at any point. Particularly, since $\displaystyle\int_{-1}^{1} \delta_S(x)dx$, we can consider it (in a sense) a derivative of $\operatorname{sign}(x)$ function. This shows us that an infinite value of the derivative at a point indicates a jump the function makes there.

Expanding $\omega^p$ as a divergent integral, Surreal Delta function can be raised to any power:

$$\delta_S(x)^p = \begin{cases} \frac{p}{\pi^p} \int_0^\infty x^{p-1}dx, & \text{if } x = 0 \\ 0, & \text{if } x \neq 0. \end{cases} \tag{4}$$

It should be noted that the Surreal Delta function is not differentiable in terms of surreal-valued functions.

## 3. Integration of a surreal value over a subset of reals

Generalizing relation (3) to any analytic function $f$ of $\omega$, we have:

$$\int_{-\infty}^{\infty} \left( \begin{cases} f(\delta_S(0)) & x = 0 \\ 0 & x \neq 0 \end{cases} \right) dx = f'(a) + \frac{1}{\pi} \int_a^\infty f''(x)\,dx, \tag{5}$$

for arbitrary real $a$ where the expression converges. This is equal to $f'(\delta_S(0))$, if $f'(\delta_S(0)) \in \mathbb{J}$ according to the definition (1).

Noticing that our H-field is equipped with operation of derivation $\partial()$, we can re-write it like this:

$$\int_{-\infty}^{\infty} \left( \begin{cases} w & x = 0 \\ 0 & x \neq 0 \end{cases} \right) dx = \pi\partial(w), \tag{6}$$

where $w \in \mathbb{J}$. But if we are integrating over a set where the expression under the integral can take surreal values at several points, it is natural to multiply this result by the numerosity of the support:

$$\int_S w dt = \pi N(S) \partial(w) \quad (7)$$

One can easily notice that surreal integration defined this way does not generally satisfy the linearity against an infinite factor $\alpha$:

$$\int_S (\alpha f) \neq \alpha \int_S f$$

For instance, if $u(x)$ is a unit impulse function $u(x)=\begin{cases} 1 & x=0 \\ 0 & x \neq 0 \end{cases}$, then

$$\int_{-1}^{1} \omega u(x) dx = \pi \neq \omega \int_{-1}^{1} u(x) dx = 0$$

In our opinion, the struggle to keep this rule intact is the reason for many failed attempts to define the surreal integration. Still, there are ongoing efforts by other authors to keep this rule [8].

As a consequence,

$$\frac{1}{\pi} \int_S \omega dt = \int_S \delta_S(0) dt = N(S). \quad (8)$$

## 4. Numerosity of the continuum

A NOTE ON THE NOTATION. Since the surreal integration may depend on whether a single point is included in the integration interval, from now on, in the context of real limits $a, b \in \mathbb{R}$, we use notation $\int_a^b$ to denote an integral where the limits of integration, $a$ and $b$ are half-included, particularly, $\int_a^b = \frac{1}{2}\left(\int_{(a,b)} + \int_{[a,b]}\right)$. The same refers to the case when only one of the limits is real.

FIRST UNCOUNTABLE ORDINAL. The automorphisms of the surreal field give us some freedom in choosing the first uncountable ordinal $\omega_1$. We have freedom to define it, using formula (8), as the numerosity of the reals on the interval $[0,1)$:

$$\omega_1 = \int_{[0,1)} \delta_S(0) dt = \frac{1}{\pi} \int_{[0,1)} \omega dt = \frac{1}{\pi} \int_0^1 \omega dt \quad (9)$$

Immediately, one can see from the formula (7) that $\omega_1$ can be represented also this way:

$$\omega_1 = \frac{1}{\pi} \int_0^{\infty} \ln \omega dt \quad (10)$$

It is important that only one of the limits of the integration to be included (or the both half-included), as we previously agreed, because the finite part of a numerosity corresponds to the Euler's characteristic of a set, and we want the finite part of $\omega_1$ to be zero.

## 5. Integrating a $\mathbb{J}+\mathbb{R}$-valued function over $\mathbb{J}+\mathbb{R}$.

EXPLICIT FORMULA. We assume that the Newton-Leibniz formula is valid for integration over $\mathbb{J}+\mathbb{R}$. As such, symbolically obtained antiderivative of an analytic function allows to calculate definite integrals on surreals. Based on this assumption and on formula (5), we propose the following explicit formula. Let $F(x)\in\mathbb{J}+\mathbb{R}$ to be a surreal-valued function with only purely-infinite and finite parts, defined on surreal domain and $a,b\in\mathbb{J}+\mathbb{R}$. Then,

$$\int_a^b F(t)dt=\pi\omega_1\int_{a(\text{formal})}^b \partial(F(t))dt+\int_a^b \operatorname{fin} F(t)dt, \tag{11}$$

where the left-hand part should be understood as integration over $t\in\mathbb{J}+\mathbb{R}$, the first integral on the right-hand side should be taken formally, the way as if all the surreal numbers under the integral were real constants, and with help of the Newton-Leibniz formula. In the second integral, the notation $\operatorname{fin}$ means taking the finite part. For the second integral also, the Newton-Leibniz formula can be used.
It should be noted that an interval of $\mathbb{J}+\mathbb{R}$ numbers is locally equivalent to an interval of reals. If one wants to integrate over a more or less dense set, one should use a more general formula (7).

MATHEMATICA CODE. Here is the Mathematica code for integration of an arbitrary surreal function over surreal domain using the formula (11). For simplicity, here `w` and `W` stand for $\omega$ and for $\omega_1$ respectively (if you want to use an arbitrary surreal in input, you should specify it as A[w] or B[w] to distinguish it from reals).

```
f[x_]:=Exp[x*w];
a=0;
b=α;
f[x_]=f[x]/. W->W[w];
Unprotect[Power];0^0=1;Protect[Power];
int=Pi*W*Integrate[D[f[t],w],{t,a,b}]+Integrate[Fin[f[t]],{t,a,b}];
int=FullSimplify[Normal[If[a==b,Pi*D[f[a],w],int,int]/. W[w]:>W/. Derivative[1]
[W][w]->∂[W]/. FullForm[Derivative[d_][W]]:>(∂^d)[W]]];
Print[Inactivate[Integrate[f[x]  /.  W[w]->W  /.  Derivative[1][W][w]->∂[W]  /.
FullForm[Derivative[d_][W][w]]:>(∂^d)[W],{x,a,b}],Integrate],"=",int]
```

It leaves the integration of the finite part unevaluated because finding the finite part of a surreal number `Fin[f[t]]` is a non-trivial problem (although straightforward if the Conway normal form of the expression is known) but for some numbers in $\mathrm{No}(\omega_1)$ one can try the regularization via Laplace transform:

```
Fin[e_] := FullSimplify[(PowerExpand[e] /. Log[w] -> 0 /. 1/w -> 0 /.
1/w^(p_) :> 0 /. w -> 0) +
             Limit[Evaluate[LaplaceTransform[D[e, w], w, x]] +
Evaluate[LaplaceTransform[D[e, w], w, -x]], x -> 0]/2 /.
        Infinity -> 0 /. -Infinity -> 0 /. FullForm[Derivative[d_][W][p_]] :> 0]
```

Alternatively, one can utilize the built-in Mathemaca’s regularization machinery:

```
Fin[e_] := (PowerExpand[f[t]] /. W[w] -> 0 /. Log[w] -> 0 /. 1/w -> 0 /.
1/w^(p_) :> 0 /. w -> 0 /. 1/W[w] -> 0) +
          Limit[FullSimplify[s*Sum[D[e, w] /. w -> s*w, {w, 1, Infinity},
Regularization -> "Dirichlet"]], s -> 0]
```

The list of some surreal indefinite integrals is given in Appendix 2.

## 6. Mapping divergent series to surreals

We have provided a mapping of some divergent integrals into surreals. Let us consider mapping of divergent series to surreals.

Let S(x) the Nørlund principal solution for the summation operator: $S(t)=\sum_{k=0}^{t-1} a_k$.

Now consider the following integral: $\int_x^{x+1} S(t)dt$. If its germ at $x\to\infty$ corresponds to a sum of a purely-infinite surreal number and a real number, then we postulate equivalence:

$$\sum_{k=0}^{\infty} a_k=\int_\omega^{\omega+1} S(t)dt \tag{12}$$

This definition keeps the regularized part of the series as the finite part of the resulting surreal number.

For instance, we can derive that $\sum_{k=0}^{\infty} k^p=\frac{\omega^{p+1}}{p+1}+\zeta(-p)+0^p$, for $p\geq 0$.
In Mathematica, this can be obtained with the following code:

```
Expand[FullSimplify[Integrate[Sum[f[k], {k, 0, x - 1}], {x, ω, ω + 1}]]]
```

where `f[k]` is a term of the series.

A list of surreal values of divergent series can be found in Appendix 3.

Consequently, we can represent infinite products as surreal numbers: $\prod_{k=1}^{\infty} f(k)=\exp\left(\sum_{k=1}^{\infty}\ln f(k)\right)$. For example,

$$\prod_{k=1}^{\infty} k=\exp\left(\sum_{k=1}^{\infty}\ln k\right)=\exp(-\omega+\omega\log(\omega)+\frac{1}{2}\log(2\pi))=\sqrt{2\pi}e^{-\omega}\omega^{\omega}$$

## 7. Comparison with prior works.

Since there are several previous works on numerosities, it is worth to highlight differences between the proposed definitions. Our aim among other things was to make a definition, in which the finite part of the numerosity of any subset of integers, particularly, the natural numbers would coincide with the regularized value (using Dirichlet or Zeta regularization) of the sum of the indicator function over the set, for instance:

$$\operatorname{fin} N(\mathbb{N}) = \operatorname{reg}\sum_{k=1}^{\infty} 1 = -1/2.$$

The numerosity of natural numbers, using our definition, corresponds to the germ of the function $f(x) = x - 1/2$ at infinity, which in turn corresponds to the surreal number $\omega - 1/2$. Similarly, definite values can be assigned to the numerosities of even or odd positive numbers: $N(odd_+) = \omega/2$ and $N(even_+) = \omega/2 - 1/2$, so their difference is ½. This result coincides with that of James Propp[7]. In comparison, Benci at al in their work[1] (p.63) concede that their definition has an arbitrariness depending on the choice of ultrafilter, and can be defined as $num(odd_+) = num(even_+)$ or $num(odd_+) = num(even_+) + 1$.

In general, other works often do not define or neglect the finite parts of the numerosities, effectively treating numerosities as growth rates.

An exception to this observation is the article by James Propp[7] on Wordpress. There he suggests a formula for finding the finite part of the numerosity of a sequence of reals using Abel's regularization:

$\operatorname{fin} N(\{a_n\}) = \lim_{q\to 1^-} \sum_{n=0}^{\infty} q^{a_n}$. Our results coincide with those of James Propp in the finite part.

Later a team (Aaron Abrams et al, including James Propp) expanded the theory [9] although this definition is not translation-invariant for bounded and even finite sets (for instance, numerosity of a singleton set can vary by an infinitesimal value).

Trlifajová[4] defines a "size" of the set, $\sigma(S)$, a measure, similar to numerosity. It is defined as the sequence of partial sums over the set's indicator function, and as such, can be a sequence itself. Operations on the "sizes" are defined as element-wise operations on sequences. In the table below for uniformity, we still will use the notation $N(S)$ for all definitions.

The following table summarizes the differences in the definitions (N/D stands for "not defined"):

| | Our definition | Propp | Trlifajová | Benci | Lynch |
|---|---|---|---|---|---|
| $\operatorname{fin} N(\mathbb{N})$ | -1/2 | -1/2 | N/D | N/D | 0 |
| $N(odd_+) - N(even_+)$ | 1/2 | 1/2 | $(1,0,1,0,1,0,...)$ | arbitrary (0 or 1) | 0 |
| Basic infinite constant used in notation | $\omega = N(\mathbb{Z})/2$ $= N(\mathbb{N}) + 1/2$ | - | $\alpha = N(\mathbb{N})$ | $\alpha = N(\mathbb{N})$ | $\omega = N(\mathbb{N})$ |

# 8. Discussion

In this work we proposed a way to map numerosities of certain subsets of reals to surreal numbers and used postulated equivalence between surreals and divergent improper integrals to interpret them in Laplace transform sense, which allowed us to derive certain integration properties of surreals. Still there remains a huge room for generalization of our principles.

For numerosities of many sequences our formula does not provide a closed-form expression, so that certain approximation techniques are to be developed. A formula for sequences with accumulation points is still to be found.

Numerosities of the sets in higher-dimensional space is another still unsolved problem. Already in $\mathbb{R}^2$, a bounded set can be made into its own proper subset by a simple rotation. This hints that a numerosity-like measure following Euclid's principle cannot be always rotation-invariant. In any case, the finite part of a numerosity of a set, corresponds to its Euler's characteristic.

The numerosities of the dense countable sets is another point of uncertainty. In her work[4] Trlifajová estimates the "size" of the unit interval of rational numbers:

$$\frac{3}{10} \cdot \alpha^2 < \sigma((0,1]_{\mathbb{Q}}) < \frac{\alpha^2 - \alpha}{2}$$

In Benci we encounter: "In particular, it seems there is no definitive way to decide whether $n_\alpha((0,1]_{\mathbb{Q}}) \geq \alpha$ or $n_\alpha((0,1]_{\mathbb{Q}}) \leq \alpha$. So, in absence of any reason to choose one of the two possibilities, we go for the simplest option $n_\alpha((0,1]_{\mathbb{Q}}) = \alpha$"[10] (Benci & Di Nasso 2019, p. 291).

On the other hand, based on the formula for the Dirichlet function $\mathbf{1}_{\mathbb{Q}}(x) = \lim_{k\to\infty}\left(\lim_{j\to\infty}(\cos(k!\pi x))^{2j}\right)$, we could represent it as a sum of lattices and naively conjecture $N(\mathbb{Q}) = 2\omega!$ and $N([0,1)_{\mathbb{Q}}) = \Gamma(\omega)$. Using similar ideas, we can naively derive that the numerosity of the dyadic rationals could be $\omega 2^\omega$.

All these ideas though, need assumption for countable additivity for numerosity and seemingly conflict with the Euclid's principle as one can construct (in all of the proposed theories) a discrete sequence of rationals with numerosity greater than $2\omega!$ and $\omega 2^\omega$ which would still be a subset of the rational numbers. Maybe, some kind of Adele ordering would be needed for finding the numerosity of the rationals.

APPENDIX 1. Numerosities of real sequences. The following is a list of calculated full numerosities (that is, with infinitesimal parts included) $N$ and refined numerosities (that is without infinitesimal parts), $\overline{N}$ of various subsets of reals ($k \in \mathbb{N}_0$):

- $N(\mathbb{N}) = N(\{k+1\}) = \omega - 1/2$
- $N(\mathbb{Z}) = 2\omega$
- $N(\{2k\}) = \dfrac{\omega}{2} + \dfrac{1}{2}$
- $N(\{2k+1\}) = \dfrac{\omega}{2}$
- $N(\{ak+b\}) = \dfrac{\omega}{a} + \dfrac{1}{2} - \dfrac{b}{a}$
- $N(\{k^2\}) = \dfrac{1}{6}\sqrt{36\omega+3} + \dfrac{1}{2};$ $\qquad \overline{N} = \sqrt{\omega} + \dfrac{1}{2}$
- $N(\{1/3 + k + k^2\}) = \sqrt{\omega}$
- $N(\{k^4\}) = \dfrac{1}{30}\sqrt{30\sqrt{900\omega+30}+225} + \dfrac{1}{2};$ $\qquad \overline{N} = \sqrt[4]{\omega} + \dfrac{1}{2}$
- $N(\{a^k\}) = \dfrac{\ln\left(\frac{(a-1)\omega}{\ln a}\right)}{\ln a}$

APPENDIX 2. Antiderivatives of some surreal functions ($c \in \mathbb{R}, c > 0$ and $\alpha \in \mathrm{No}, \alpha \neq 0$).

- $\displaystyle\int_0^\alpha \omega^c dt = \pi\alpha c\omega_1\omega^{c-1}$
- $\displaystyle\int_0^\alpha \omega^x dx = \frac{\pi\left((\alpha\ln\omega - 1)\omega^\alpha + 1\right)\omega_1}{\omega(\ln\omega)^2}$
- $\displaystyle\int_0^\alpha x^\omega dx = \frac{\pi\alpha^{\omega+1}((\omega+1)\ln\alpha - 1)\omega_1}{(\omega+1)^2}$
- $\displaystyle\int_0^\alpha e^\omega dx = e^\omega\pi\alpha\omega_1$
- $\displaystyle\int_0^\alpha c^\omega dx = c^\omega\pi\alpha\omega_1\ln c$
- $\displaystyle\int_0^\alpha e^{cx\omega} dx = \frac{\pi\left(e^{c\alpha\omega}(c\alpha\omega - 1) + 1\right)\omega_1}{c\omega^2}$
- $\displaystyle\int_0^\alpha \ln\omega dx = \frac{\pi\alpha\omega_1}{\omega}$
- $\displaystyle\int_0^\alpha \omega^\omega dx = \alpha\omega^\omega\omega_1\pi(\ln\omega + 1)$
- $\displaystyle\int_0^\alpha \omega^{\omega_1} dx = \alpha\pi\omega_1\omega^{\omega_1-1}\left(\omega\partial\left(\omega_1\right)\ln\omega + \omega_1\right)$
- $\displaystyle\int_0^\alpha \omega_1^\omega dx = \alpha\pi\omega_1^\omega\left(\omega\partial\left(\omega_1\right) + \omega_1\ln\omega_1\right)$

APPENDIX 3. Surreal values of some divergent series ($p \in \mathbb{R}, p \geq 0$).

- $\sum_{k=0}^{\infty} k^p = \frac{\omega^{p+1}}{p+1} + \zeta(-p) + 0^p$
- $\sum_{k=0}^{\infty} p^k = \frac{p^\omega}{\ln p} + \frac{1}{1-p}$
- $\sum_{k=1}^{\infty} \ln k = \omega \ln \omega - \omega + \frac{1}{2} \ln(2\pi)$
- $\sum_{k=1}^{\infty} \frac{1}{k} = \ln \omega + \gamma$
- $\sum_{k=1}^{\infty} \psi(k) = \ln \Gamma(\omega) + \frac{1}{2} - \frac{1}{2} \ln(2\pi)$

Here $\psi(x)$ is the digamma function and $\gamma$ is the Euler-Mascheroni constant.